\documentclass[AMA,STIX1COL]{WileyNJD-v2}
\usepackage{moreverb}

\newcommand\BibTeX{{\rmfamily B\kern-.05em \textsc{i\kern-.025em b}\kern-.08em
T\kern-.1667em\lower.7ex\hbox{E}\kern-.125emX}}

\articletype{Article Type}%

\received{<day> <Month>, <year>}
\revised{<day> <Month>, <year>}
\accepted{<day> <Month>, <year>}

\begin{document}

\title{Non-Hydrostatic Model for Simulating Moving Bottom-Generated Waves: A Shallow Water Extension with Quadratic Vertical Pressure Profile}

\author[1,2]{Kemal Firdaus*}

\author[1,2]{J\"orn Behrens}

\authormark{FIRDAUS and BEHRENS}

\address[1]{Department of Mathematics, Universit\"at Hamburg, Bundesstrasse 53-55, 20146 Hamburg, Germany}

\address[2]{Center for Earth System Research and Sustainability (CEN), Universit\"at Hamburg, Bundesstrasse 53-55, 20146 Hamburg, Germany}

\corres{*Kemal Firdaus, Department of Mathematics, Universit\"at Hamburg, Bundesstrasse 53-55, 20146 Hamburg, Germany.\\ \email{kemal.firdaus@uni-hamburg.de}}


\abstract[Abstract]{We formulate a depth-averaged non-hydrostatic model to solve wave equations with generation by a moving bottom. This model is built upon the shallow water equations, which are widely used in tsunami wave modelling. An extension leads to two additional unknowns to be solved: vertical momentum and non-hydrostatic pressure. We show that a linear vertical velocity assumption turns out to give us a quadratic pressure relation, which is equivalent to Boussinesq-type equations. However, this extension involves a time derivative of an unknown parameter, rendering the solution by a projection method ambiguous. In this study, we derive an alternative form of the elliptic system of equations to avoid such ambiguity. The new set of equations satisfies the desired solubility property, while also consistently representing the non-flat moving topography wave generation. Validations are performed using several test cases based on previous experiments and a high-fidelity simulation. First, we show the efficiency of our model in solving a vertical movement, which represents an undersea earthquake-generated tsunami. Following that, we demonstrate the accuracy of the model for landslide-generated waves. Finally, we compare the performance of our novel set of equations with the linear and simplified quadratic pressure profiles.} 


\maketitle

\section{Introduction}\label{sec1}
Large-scale ocean dynamics are affected by long-period gravity waves. Numerous studies have been conducted to model such waves. One of the most popular models is given by the shallow water equations (SWE) due to its simplicity and capacity to explain a wide range of wave dynamics. A prominent example of the use of SWE is in tsunami wave propagation \cite{behrens2015, harig2008, FIRDAUS2022, MAGDALENA2022}. However, the model is limited to situations being in hydrostatic balance, where the pressure is assumed to be hydrostatic only. Non-hydrostatic forces, however, can be necessary in some situations, particularly when dealing with landslides and slow earthquake-generated tsunamis \cite{glimsdal2013}. In fact, submarine landslide contributes to $24\%$ of all fatal tsunamis that exist over the last century \cite{reid2023}. For such cases, the fluid depth and wavelength ratio become higher, giving a dispersion effect, which can not be handled with a model limited by the hydrostatic assumption. In that case, Boussinesq-type equations have been successfully applied \cite{GOBBI1999, FUHRMAN2009747, PANDA2014}.

Boussinesq-type equations are derived by applying an asymptotic expansion on the non-dimensional velocity potential, which is based on the incompressible Euler equations, as a power series in the vertical coordinate. Previous studies used different expansion orders and formulations, giving various types of Boussinesq equations \cite{peregrine1967, nwogu1993, madsen1998, madsen2003}. Based on their expansion order, the models are described as weakly, moderately, or fully dispersive and nonlinear. A review of the development of Boussinesq-type equations is given in \cite{brocchini2013}. Despite their physical capabilities, higher-order derivatives and mixed space-time derivatives appear in such sets of equations, which may lead to difficulties in their numerical discretisation. Therefore, computationally advantageous models are still under investigation.

Considering the SWE again, a derivation can be achieved by applying depth integration to the Euler equations under the assumption of hydrostatic pressure. Then, a non-hydrostatic extension is constructed by decomposing the pressure into a hydrostatic and a non-hydrostatic component \cite{casulli1998, stansby1998}. This yields two additional quantities: vertical momentum and non-hydrostatic pressure. To close the system of equations, two additional equations are required, which are the vertical momentum and divergence constraint equations. In the end, a vertical pressure relation is required to achieve a purely averaged system. Some studies assumed a linear vertical pressure relation \cite{walters2005, stelling2003, yamazaki2009}, while also assuming a linear vertical velocity. Alternatively, \cite{seabra-santos1987, Jeschke2017} derived a quadratic pressure relation based on the linear velocity assumption. The latter relation is proven to be consistent with the Boussinesq-type equations and has also been used for previous studies on static bottom cases \cite{DEMPWOLFF2024, wang2020}. This extension is then solved with a projection-based method, which was initially introduced by \cite{chorin1968} to numerically solve the incompressible Navier-Stokes equations. The pressure decomposition method allows us to solve a non-hydrostatic correction based on the result of the hydrostatic SWE solver, which is therefore called a predictor. Following the SWE predictor, we correct the solution by solving the non-hydrostatic pressure and vertical momentum upon the calculated predictor. At this point, we will end up solving an elliptic system of equations in each time step. However, as described in \cite{jeschke2018}, the quadratic pressure relation involves a time derivative of the horizontal velocity, which is unknown. 

The solution of the time derivative term is not straightforward, which is why in \cite{jeschke2018} it was argued that the term could be neglected. In this study, we propose an alternative formulation of the pressure relation avoiding the time derivative. 

To solve the derived model, we discretize it with a discontinuous Galerkin (DG) method in space. This model satisfies desired properties in geophysical flow modelling, such as conservation, mesh flexibility, and accuracy\cite{bernard2007, giraldo2008, DAWSON2013}. The locality of its shape function also allows for parallelization and application of adaptivity \cite{VATER20151, MULLER2013}. For the time-stepping, we use a second-order Runge-Kutta scheme for solving the predictor and an implicit backward Euler for the correction. Moreover, to solve the elliptic system of equations, we use the local discontinuous Galerkin (LDG) method. This method is initially an extension of the RKDG method for a time-dependent non-linear advection-diffusion system developed by \cite{cockburn1998}. Here the LDG is used to solve the elliptic equation, similar to the Poisson equation solvers in \cite{cockburn2003, CASTILLO20061307}. More general stable numerical fluxes for three types of boundary conditions are derived in \cite{jeschke2018}. However, the analysis was restricted to flat bottom cases, to achieve the required structure properties. This work will also show that taking care of the topography steepness and movement will preserve the desired properties. Hence, the derived flux can be adapted directly.

The rest of the paper is structured as follows. The non-hydrostatic extension for the SWE for moving bottom-generated waves will be derived in Section \ref{sec2}. Along with that, we also provide the alternative form of the pressure relation to avoid the time-derivative of the unknown in the non-hydrostatic pressure term. Following that, Section \ref{sec3} will discuss further the numerical method to solve the derived model. There, we will show the elliptic system of equations that need to be solved, which fulfils the desired properties for the flux. Finally, we then compare our simulation results with two experimental data and a former simulation in Section \ref{sec4} and conclude it in Section \ref{sec5}.

\section{Mathematical Model}\label{sec2}
In this section, we derive the dispersive extension for the SWE along with its non-hydrostatic pressure relation in a one-dimensional case. First of all, we provide the depth-averaged model based on the incompressible Euler equations, where we also include the bottom movement term. Following that, we come up with the alternative form of bottom non-hydrostatic pressure relation with the averaged non-hydrostatic pressure. Again, the latter part is crucial to avoid confusion in applying the projection method in solving the elliptic system of equations.

\subsection{The non-hydrostatic extension for SWE with a moving bottom}
\begin{figure}[h]
    \centering
    \includegraphics[width=0.4\linewidth]{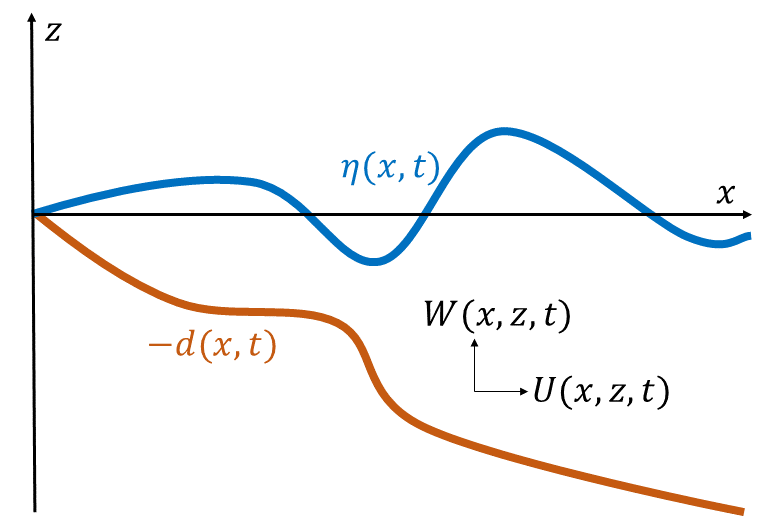}
    \caption{The illustration of the model.}
    \label{fig:swe}
\end{figure}

We begin with the two-dimensional incompressible Euler equations of motion, which can be written as
\begin{equation}\label{eq:euler_1}
    U_x+W_z= 0,
\end{equation}
\begin{equation}\label{eq:euler_2}                          
    U_t+(U^2)_x+(UW)_z=-\frac{1}{\rho}P_x,
\end{equation}
\begin{equation}\label{eq:euler_3}                          
    W_t+(UW)_x+(W^2)_z=-\frac{1}{\rho}P_z-g,
\end{equation}
with $(U,W)^T$ being the velocity in $x,z$-direction respectively, $\rho$ is the fluid density, $P$ is the total pressure, and $g$ is the gravitational acceleration. For an illustration of the involved variables here, see Fig. \ref{fig:swe}. Note that Eq. \ref{eq:euler_1} corresponds to the continuity equation (conservation of mass), while Eq. \ref{eq:euler_2} and \ref{eq:euler_3} correspond to the horizontal and vertical momentum balance equations respectively. To be solvable, these equations are closed by the kinematic boundary conditions at the surface and fluid bottom and the pressure condition, which respectively can be written as
\begin{equation}\label{eq:kinematic_surface}
    W_\eta:=W(x,\eta,t) = \eta_t+U|_{z=\eta}\eta_x,
\end{equation}
\begin{equation}\label{eq:kinematic_bottom}
    W_{-d}:=W(x,-d,t) = -d_t-U|_{z=-d}d_x,
\end{equation}
\begin{equation}\label{eq:pressure_surface}
    P|_{z=\eta}=0.
\end{equation}
Note that we include the bottom movement by a non-vanishing term $d_t$ in the bottom kinematic boundary condition. This term is often neglected, especially when either instantaneous bottom uplift or surface perturbation is assumed as a source. In ordinary SWE, $P$ is purely hydrostatic, while here, we assume that it has a non-hydrostatic component called $P^{nh}$. So, in addition to the hydrostatic pressure $P^{hy}$, we write
\begin{equation}\label{eq:pressure_expand}
    P = P^{hy}+P^{nh} = \rho g(\eta-z)+P^{nh},
\end{equation}
which implies $P|_{z=\eta} = P^{hy}|_{z=\eta} = P^{nh}|_{z=\eta}=0$. 

Integrating Eq. \ref{eq:euler_1}, averaging it over the fluid depth $h=\eta+d$, and making use of the kinematic boundary conditions (\ref{eq:kinematic_surface}), (\ref{eq:kinematic_bottom}) we obtain
\begin{equation}\label{eq:swe_mass}
   h_t+(hu)_x=0,
\end{equation}
where we use the depth-averaged unknown horizontal velocity, defined as $u=\frac{1}{h}\int_{-d}^\eta Udz$. In a similar manner, and involving the pressure relations (\ref{eq:pressure_surface}), (\ref{eq:pressure_expand}) we can deduce from Eq. \ref{eq:euler_2} and \ref{eq:euler_3},
\begin{equation}\label{eq:swe_momentum_1}
    (hu)_t+\left(hu^2+\frac{ g}{2}h^2\right)_x=ghd_x-\frac{1}{\rho}(hp^{nh})_x+\frac{1}{\rho}P^{nh}|_{z=-d}d_x
\end{equation}
\begin{equation}\label{eq:swe_momentum_2}
    (hw)_t+(huw)_x=\frac{1}{\rho}P^{nh}|_{z=-d},
\end{equation}
with $w=\frac{1}{h}\int_{-d}^\eta Wdz,\quad p^{nh}:=\frac{1}{h}\int_{-d}^\eta P^{nh}dz$. To achieve such depth-averaged nonlinear terms, we assume a small vertical variation in the horizontal velocity, which holds for shallow water cases. With the same assumption, we can approximate $\boldsymbol{U}$ with $\boldsymbol{u}$ in the kinematic boundary conditions, that we achieve the depth-averaged form of the kinematic boundary conditions
\begin{equation}\label{eq:kinematic_surface_app}
    W|_{z=\eta} = \eta_t+u\eta_x,
\end{equation}
\begin{equation}\label{eq:kinematic_bottom_app}
    W|_{z=-d} = -d_t-ud_x.
\end{equation}

Following \cite{Jeschke2017}, we transform the boundary conditions (\ref{eq:kinematic_surface})  and (\ref{eq:kinematic_bottom}) to achieve purely depth-averaged equations. Assume the vertical velocity $W$ to be linear, we obtain the depth average as $w = \frac{1}{2}(W|_{z=\eta}+W|_{z=-d})$, which leads to
\begin{align*}
    w &= \frac{1}{2}(W|_{z=\eta}+W|_{z=-d})=\frac{1}{2}(W|_{z=\eta}-W|_{z=-d}+2W|_{z=-d})\\
    &=\frac{1}{2}(h_t+uh_x+2W|_{z=-d}) = \frac{1}{2}(-hu_x+2W|_{z=-d})
\end{align*}
\begin{equation}\label{eq:swe_constraint}
    \implies2hw +hu(2d-h)_x+2hd_t = -h(hu)_x
\end{equation}
At this point, we have formulated the non-hydrostatic extension for SWE as in Eq. \ref{eq:swe_mass}, \ref{eq:swe_momentum_1}, \ref{eq:swe_momentum_2}, and \ref{eq:swe_constraint}.
This set of equations has already been shown to be equivalent to Boussinesq-type equations in \cite{Escalante2020}.

\subsection{The non-hydrostatic pressure profile}
As can be seen in the previous section, we are still required to write the non-hydrostatic pressure at the bottom with its depth-averaged value to achieve an averaged system. Some studies assume a linear pressure profile, such that 
\begin{equation}\label{eq:linear_pressure}
    P^{nh}|_{z=-d}=2p^{nh}.
\end{equation}
On the other hand, if we stick to the linear vertical velocity profile, which we used in order to derive the divergence constraint, we obtain the quadratic vertical pressure relation \cite{seabra-santos1987, Jeschke2017}. For static bottom cases, it holds
\begin{equation}
    P^{nh}|_{z=-d} = \frac{3}{2}p^{nh}-\frac{\rho h}{4}(d_x(u_t+uu_x)+u^2d_{xx}).
\end{equation}
However, as stated in \cite{Jeschke2017}, it is impossible to solve such a quadratic relation with a projection method. A practical way around is to neglect some terms and thereby avoid the time-dependent unknown values, obtaining the desired flux, which we end up with
\begin{equation}\label{eq:quadratic_pressure_old}
    P^{nh}|_{z=-d} = \frac{3}{2}p^{nh}.
\end{equation}
This approach can be found in \cite{Jeschke2017, jeschke2018, ESCALANTE2019}. Here, we re-investigate the quadratic pressure relation, which we will apply to moving bottom cases. 

First of all, recall the assumption of linear vertical velocity, such that
\begin{equation}
    W(z) = W|_{z=-d}+\frac{1}{h}(W|_{z=\eta}-W|_{z=-d})(z+d).
\end{equation}
Substituting \ref{eq:kinematic_surface_app} and \ref{eq:kinematic_bottom_app}, we obtain
\begin{equation}\label{eq:wt}
    W_t = -d_{tt}-u_td_x-ud_{xt}-u_{xt}(z+d)-u_xd_t,
\end{equation}
\begin{equation}\label{eq:uwx}
    uW_x = -ud_{xt}-uu_xd_x-u^2d_{xx}-uu_{xx}(z+d)-uu_xd_x,
\end{equation}
\begin{equation}\label{uxw}
    -u_xW=u_xd_t+uu_xd_x+u_x^2(z+d).
\end{equation}
Substituting these last three equations into Euler's vertical momentum equation and approximating the horizontal velocity with its depth-integrated value, we obtain
\begin{equation}
    -\frac{1}{\rho}P^{nh}_z=(z+d)(-u_{xt}-uu_{xx}+u_x^2)-d_x(u_t+uu_x)-d_{tt}-2ud_{xt}-u^2d_{xx}.
\end{equation}
Integrating this equation to evaluate the depth-averaged value and comparing it with the pressure at the bottom, gives
\begin{equation}\label{eq:pressure_init}
    P^{nh}|_{z=-d} = \frac{3}{2}p^{nh}-\frac{\rho h}{4}(d_x(u_t+uu_x)+d_{tt}+2ud_{xt}+u^2d_{xx}).
\end{equation}

Note that a time derivative of the unknown $u$ still occurs. Instead of neglecting this term, here we derive an alternative form of this pressure relation. To achieve that, note that the horizontal momentum balance can be rewritten in the non-conservative form as
\begin{equation}
    u_t+uu_x=-g\eta_x-\frac{1}{\rho h}((hp^{nh})_x-P_{-d}^{nh}d_x).
\end{equation}
By making use of the last equation, Eq. \ref{eq:pressure_init}, can alternatively be written as
\begin{equation}\label{eq:pressure_rev}
    P_{-d}^{nh}=\frac{6}{4+d_x^2}p^{nh}+\frac{d_x}{4+d_x^2}(hp^{nh})_x+\phi,
\end{equation}
where $\phi = \frac{\rho h}{4+dx^2}\left(gd_x\eta_x-d_{tt}-2ud_{xt}-u^2d_{xx}\right)$. This new form does not contain any time derivative of the unknown values. Furthermore, higher derivative order only applies to the bottom function, which is given. Therefore, this proposed relation can be directly solved with the projection method, while maintaining its equivalence with Boussinesq-type equations by not losing any terms.

\section{Numerical Method}\label{sec3}
We solve the above model numerically with a projection method, where a discontinuous Galerkin discretization is used in space. The brief idea of this method is as follows: first of all, the (simple) hydrostatic SWE are solved as predictors, then the non-hydrostatic pressure is computed based on the predicted values, finally, the predictor is corrected using the non-hydrostatic pressure values. The non-hydrostatic pressure correction terms are not solved explicitly in this section but will be discussed further down.

\subsection{Predictor step}
To obtain the predictor, we neglect the non-hydrostatic terms in our model and solve it with a second-order Runge-Kutta discontinuous Galerkin (RKDG2) method using a limiter as presented in \cite{VATER20151}. This limiter allows for wetting and drying and plays a crucial role in applying the model to landslide tsunamis, which involve inundation. Rewriting the hydrostatic form of \ref{eq:swe_mass}-\ref{eq:swe_momentum_2} in the compact conservative form, we have
\begin{equation}\label{eq:swe_compact}
    \boldsymbol{q}_t+\boldsymbol{f}(\boldsymbol{q})_x= \boldsymbol{s}(\boldsymbol{q}),
\end{equation}
with $\boldsymbol{q}: = (h,hu,hw)^T$ denote the unknowns, while $\boldsymbol{f}(\boldsymbol{q})=\left(hu,hu^2+\frac{g}{2}h^2,huw\right)^T$ and $\boldsymbol{s}(\boldsymbol{q})=(0, ghd_x, 0)^T$ denote the flux and source functions, respectively. Later on, we solve the governing equation on the discretized domain $\Omega=[x_l,x_r]$, which consists of $N$ uniform intervals $I_i=[x_i,x_{i+1}]$. Multiplying \ref{eq:swe_compact} with a test function $\phi$ and integrating it on an interval leads to the weak DG formulation
\begin{equation}
    \int_{I_i}\phi\boldsymbol{q}_tdx-\int_{I_i}\phi_x\boldsymbol{f}(\boldsymbol{q})dx+\left[\phi\boldsymbol{f}^*(\boldsymbol{q})\right]_{x_i}^{x_{i+1}}=\int_{I_i}\phi\boldsymbol{s}(\boldsymbol{q})dx.
\end{equation}
The interface flux $\boldsymbol{f}^*$ here is not generally defined due to the discontinuity, which may lead to different values on neighbouring elements. The solution of this flux is then approximated with a Riemann solver. In this case, we use the Rusanov solver (for more details and other Riemann solvers, check \cite{toro2013riemann}). A piecewise polynomial ansatz is used to discretize the solution and test the function spatially. This semi-discrete system is then solved by a total-variation diminishing (TVD) two-stage Runge–Kutta scheme, also known as Heun's method \cite{SHU1988439}. The calculated predictors here are then stored as $\tilde{\boldsymbol{q}}^{n+1}=(\tilde{h}^{n+1},\tilde{hu}^{n+1},\tilde{hw}^{n+1})^T$.

\subsection{Corrector step}
Following this predictor, an implicit Euler method is used to correct the momentum equations. On the other hand, since there is no non-hydrostatic pressure term in the mass conservative equation, the fluid depth value is adopted directly from the predictor, i.e. $h^{n+1} = \tilde{h}^{n+1}$. The momentum values are corrected based on
\begin{equation}\label{eq:correction1}
    \frac{(hu)^{n+1}-(\tilde{hu})^{n+1}}{\Delta t} = -\frac{1}{\rho}(\tilde{h}p^{nh})_x^{n+1}+\frac{d_x^{n+1}}{\rho}\left(\frac{6}{4+(d_x^{n+1})^2}(p^{nh})^{n+1}+\frac{d_x^{n+1}}{4+(d_x^{n+1})^2}(\tilde{h}p^{nh})_x^{n+1}+\phi^{n+1}\right),
\end{equation}
\begin{equation}\label{eq:correction2}
\frac{(hw)^{n+1}-(\tilde{hw})^{n+1}}{\Delta t} = \frac{1}{\rho}\left(\frac{6}{4+(d_x^{n+1})^2}(p^{nh})^{n+1}+\frac{d_x^{n+1}}{4+(d_x^{n+1})^2}(\tilde{h}p^{nh})_x^{n+1}+\phi^{n+1}\right),
\end{equation}
where $\phi^{n+1} := \phi((d,\tilde{h},\tilde{hu})^{n+1})$. 

\subsection{Elliptic pressure equation}
Dividing equation(\ref{eq:correction1}) by $\frac{4\tilde{h}^{n+1}}{\rho(4+(d_x^{n+1})^2)}$ leads to
\begin{equation}\label{eq:new_elliptic_1}
    (p^{nh})_x^{n+1}+\frac{2\tilde{h}_x^{n+1}-3d_x^{n+1}}{2\tilde{h}^{n+1}}(p^{nh})^{n+1}+\frac{\rho(4+(d_x^{n+1})^2)}{4\Delta t\tilde{h}^{n+1}}(hu)^{n+1} = \frac{4+(d_x^{n+1})^2}{4\tilde{h}^{n+1}}\left(\phi^{n+1} d_x^{n+1}+\frac{\rho}{\Delta t}(\tilde{hu})^{n+1}\right).
\end{equation}
Furthermore, making use of the discretized divergence constraint
\begin{equation}
\begin{split}
    2\left((\tilde{hw})^{n+1}+\frac{\Delta t}{\rho}\left(\frac{6}{4+(d_x^{n+1})^2}(p^{nh})^{n+1}+\frac{d_x^{n+1}}{4+(d_x^{n+1})^2}(\tilde{h}p^{nh})_x^{n+1}+\phi\right)\right)&\\
    +(2d-\tilde{h})_x^{n+1}(hu)^{n+1}+2(\tilde{h}d_t)^{n+1}&=-\tilde{h}^{n+1}(hu)_x^{n+1},
\end{split}
\end{equation}
the latter equation can be transformed into
\begin{equation}
\label{eq:new_elliptic_2}
    (hu)_x^{n+1}+\frac{3d_x^{n+1}-2\tilde{h}_x^{n+1}}{2\tilde{h}^{n+1}}(hu)^{n+1}+\frac{3\Delta t}{\rho \tilde{h}^{n+1}}(p^{nh})^{n+1}=-2d_t^{n+1}-\frac{2(\tilde{hw})^{n+1}}{\tilde{h}^{n+1}}-\frac{d_x^{n+1}(\tilde{hu})^{n+1}}{2\tilde{h}^{n+1}}-\frac{\Delta t(4+(d_x^{n+1})^2)}{2\rho\tilde{h}^{n+1}}\phi^{n+1}.
\end{equation}
It can be observed that (\ref{eq:new_elliptic_1}) and (\ref{eq:new_elliptic_2}) form an elliptic system of equations of $(p^{nh})^{n+1}$ and $(hu)^{n+1}$. Furthermore, the correction for the vertical momentum can be solved straightforwardly with the corrected pressure value. 

Following \cite{jeschke2018}, a general elliptic system of equations can be written as
\begin{equation}
    p_x+g_1p+h_1u=f_1\quad,\text{on } \Omega,
\end{equation}
\begin{equation}
    u_x+h_2p+g_2u=f_2\quad,\text{on } \Omega.
\end{equation}
In the previously mentioned study, properties of $g_1+g_2 = 0$ and $h_1,h_2>0$ are required in the derivation to achieve a suitable numerical flux. In their study, a flat bottom assumption was assumed, to achieve such properties. Our model definitely fulfils the former property, while the latter is always fulfilled, thanks to the positivity preserving scheme adopted from \cite{VATER20151}. Therefore, our novel proposed system matches the desired structure while considering the bottom divergence and time derivation. Therefore, we may apply the same numerical flux, without losing the effect of the non-flat moving bottom.

\section{Results and Discussion}\label{sec4}
To validate our model, four test cases will be presented in this section. To begin with, we simulate a wave generated by a vertically moving bar to check the ability to simulate wave propagation over a non-flat bottom. Following that, we simulate three cases of moving-bottom-generated waves. The first two are cases where the bottom movement appears over a flat bottom. These test cases are intended for the initial evaluation of a moving bottom, including one with vertical and horizontal movement. The last benchmark comprises a moving perturbation on a sloping beach, which involves inundation.

\subsection{Impulsive vertical thrust}
Here, we investigate our model's ability to simulate a wave generated by vertical bottom movement. We simulate an exponential bottom movement based on an experiment conducted by \cite{Hammack_1973} in a wave tank with a size of $31.6~m \times 0.394~m \times 0.61~m$. The tank was filled with undisturbed $h_0=0.05~m$ deep water as the initial condition. To generate waves, a moving plate with a width of $b = 0.61~m$ was placed on the upstream side, where $x=0~m$ was the upstream end being a fully reflecting wall. On the other hand, a wave dissipation system was installed at the downstream end so that no wave was reflected. Two types of movement were considered here: uplift and down draft, where the elevations were recorded at $\frac{x-b}{h_0}={0, 20, 180, 400}$ (corresponding to $x = \{0.61, 1.61, 9.61, 20.61\}~m$). The exponential bottom movement followed
\begin{equation}
    d(x,t) = h_0-\zeta_0(1-e^{-\alpha t})\mathcal{H}(b^2-x^2),
\end{equation}
where $\mathcal{H}$ is the Heaviside function and $|\zeta_0|=0.005~m$ is the maximum plate displacement. Controlling the movement velocity, $\alpha=\frac{1.11}{t_c}$ was used, where $\frac{t_c\sqrt{gh_0}}{b}=0.148$ and $0.093$ for upward and downward case respectively.

\begin{figure}[h]
    \centering
    \includegraphics[width=1\linewidth]{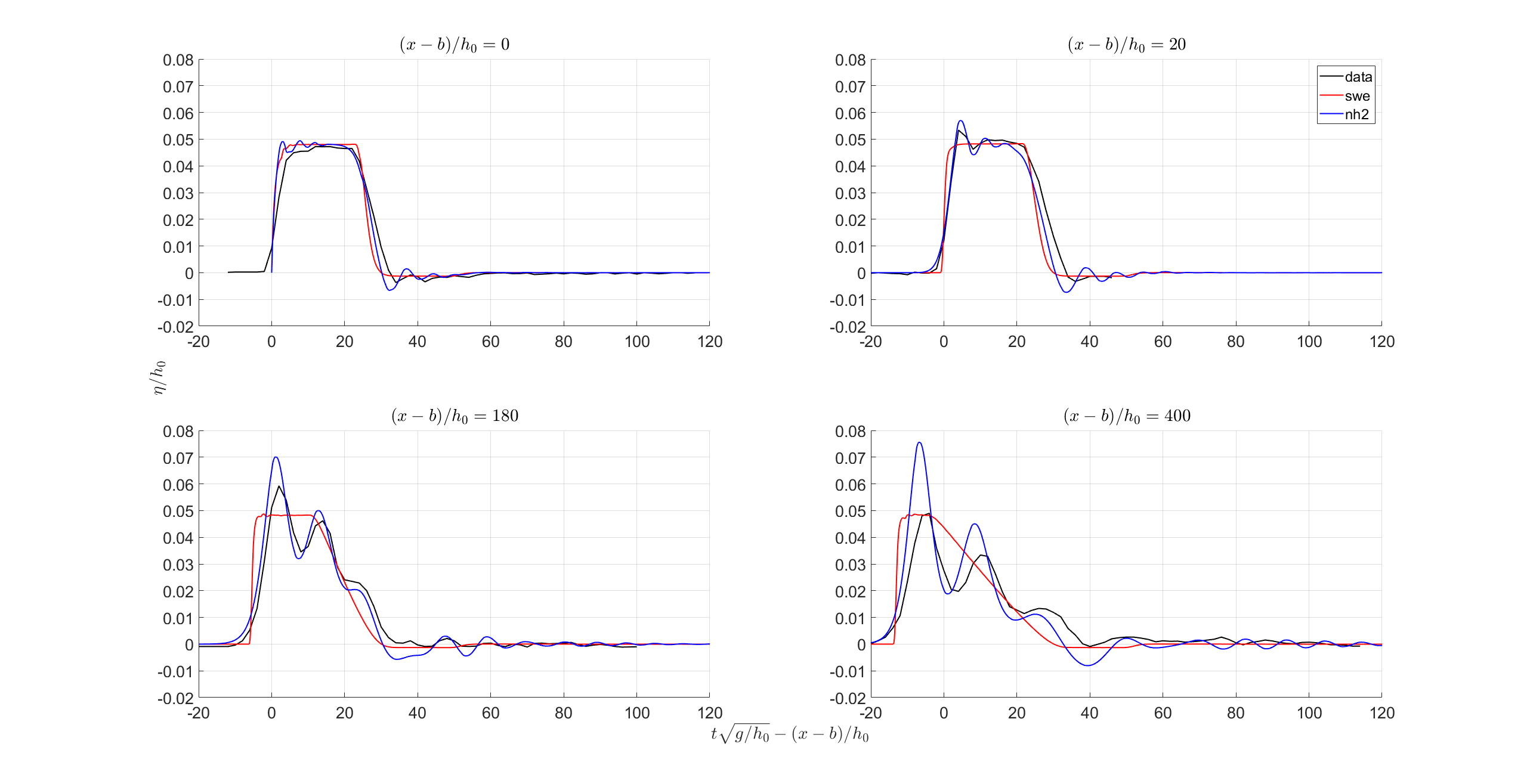}
    \caption{Comparison of the measured laboratory data (black) elevation with the simulation using hydrostatic SWE (red) and simulation using non-hydrostatic SWE (blue) for the upthrust case.}
    \label{fig:hammack_up}
\end{figure}

\begin{figure}[h]
    \centering
    \includegraphics[width=1\linewidth]{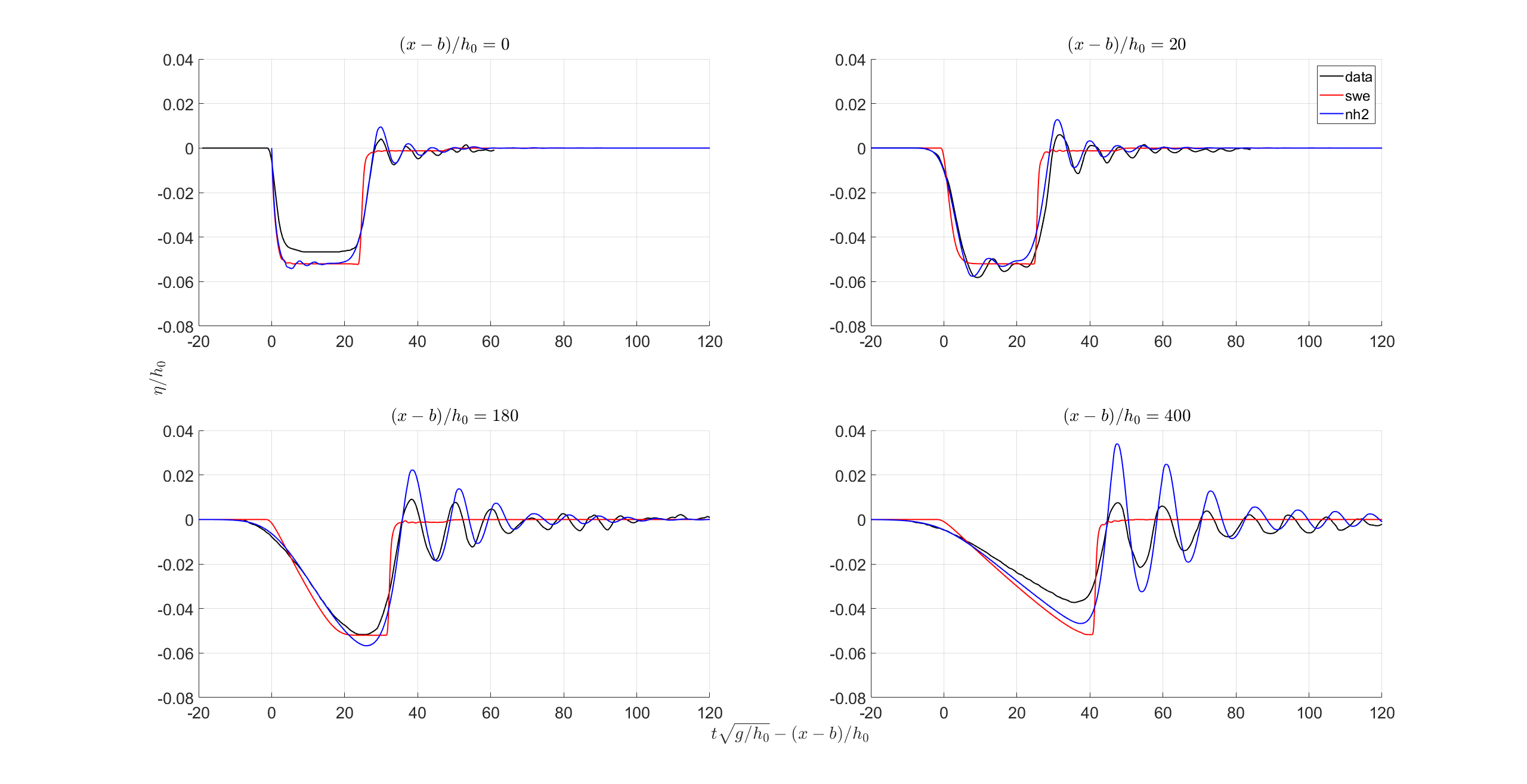}
    \caption{Comparison of the measured laboratory data (black) elevation with the simulation using hydrostatic SWE (red) and simulation using non-hydrostatic SWE (blue) for the down-thrust case.}
    \label{fig:hammack_down}
\end{figure}

For our computational setup, a time step of $\Delta t = 0.001~s$ and a spatial step size of $\Delta x = 0.025~m$ are chosen, to achieve the desired convergence. Figures \ref{fig:hammack_up} and \ref{fig:hammack_down} depict the comparison of the measured laboratory data with our simulation results for the uplift and downdraft cases respectively. Additionally, we also include the result using hydrostatic SWE. Overall, our results give a good match with the measured data. The dispersive effect also plays a significant role in both cases.

For the upward movement, our model can closely approximate the data near the generation area ($x=0.61, 1.61~m$). The hydrostatic SWE, on the other hand, still manage to give a pretty accurate waveform, even though the dispersive effect is not captured. Moreover, away from the plate, abrupt and sharp waveforms can be observed from hydrostatic SWE simulations. In contrast to that, our non-hydrostatic model captures the same waveform and period, where a shift to the right occurred in previous studies with more advanced models \cite{FUHRMAN2009747, Hammack_1973, xin024}. The excess in amplitude, which was also observed in the mentioned previous studies, might be attributed to viscous energy losses and boundary stresses in the experiments (as suggested in \cite{Hammack_1973}), which is neglected in all mentioned modelling approaches, yet leads to damping in the physical model. 

Moving to the downdraft case, we can also see that our model is able to match the data, especially near the generation area. In this case, the dispersive effect is even more notable and well captured by our model, in particular in representing the phase. As in the uplift case, the amplitude of the model is larger than in the experimental data. The hydrostatic SWE simulation fails to depict the dispersive wave effect. These two cases demonstrate the ability of the new non-hydrostatic approximation to accurately represent dispersive wave behaviour from simple bottom uplift and downdraft.

\subsection{Sliding bump over flat bottom}
In our next test case, a horizontal bottom movement is used to mimic landslide-generated waves, however, with the simplification of a flat bottom. We consider a sliding semi-elliptical block along a horizontal boundary, which is based on an experiment by \cite{whittaker2015}. This experiment was done in a $14.66~m\times0.25~m\times 0.5~m$ flume located at the Fluid Mechanics Laboratory of the Department of Civil and Natural Resources Engineering at the University of Canterbury, New Zealand. The tank was filled with a $h_0=0.175~m$ depth water, where a semi-elliptical aluminium plate, with a dimension of $0.5~m$ long, $0.25~m$ wide and $0.026~m$ thick, was placed in the middle of the wave tank ($x=0~m$). To measure the elevation, the water was also mixed with fluorescent and captured with a laser-induced fluorescent (LIF). For our comparison here, we take three cases, with parameters given in Table \ref{tab:whittaker} (The run ID also follows from the mentioned article). From a physical point of view, the landslide moves with a constant acceleration $a_0$ until achieving a velocity of $u_t$ after $t_1$, followed by a constant movement until $t_2$, and finally experiences a deceleration of $-a_0$ until it is fully stopped at $t_3$. 
\begin{table}[h!]
    \centering
    \begin{tabular}{c c c c c c c}
        \hline
         Run & $Fr$ & $a_0~(m/s^2)$ & $u_t~(m/s)$ & $t_1~(s)$ & $t_2~(s)$ & $t_3~(s)$ \\
         \hline
         6& 0.125 & $1.500$ & $0.163$ & $0.109$ & $2.109$ & $2.218$  \\
         12& 0.25 & $1.500$ & $0.327$ & $0.218$ & $2.218$ & $2.436$ \\
         18& 0.375 & $1.500$ & $0.491$ & $0.327$ & $2.327$ & $2.654$  \\
         \hline
    \end{tabular}
    \caption{Dimensional parameters of the experimental setup by \cite{whittaker2015}.}
    \label{tab:whittaker}
\end{table}

\begin{figure}[h]
    \centering
    \includegraphics[width=1 \linewidth]{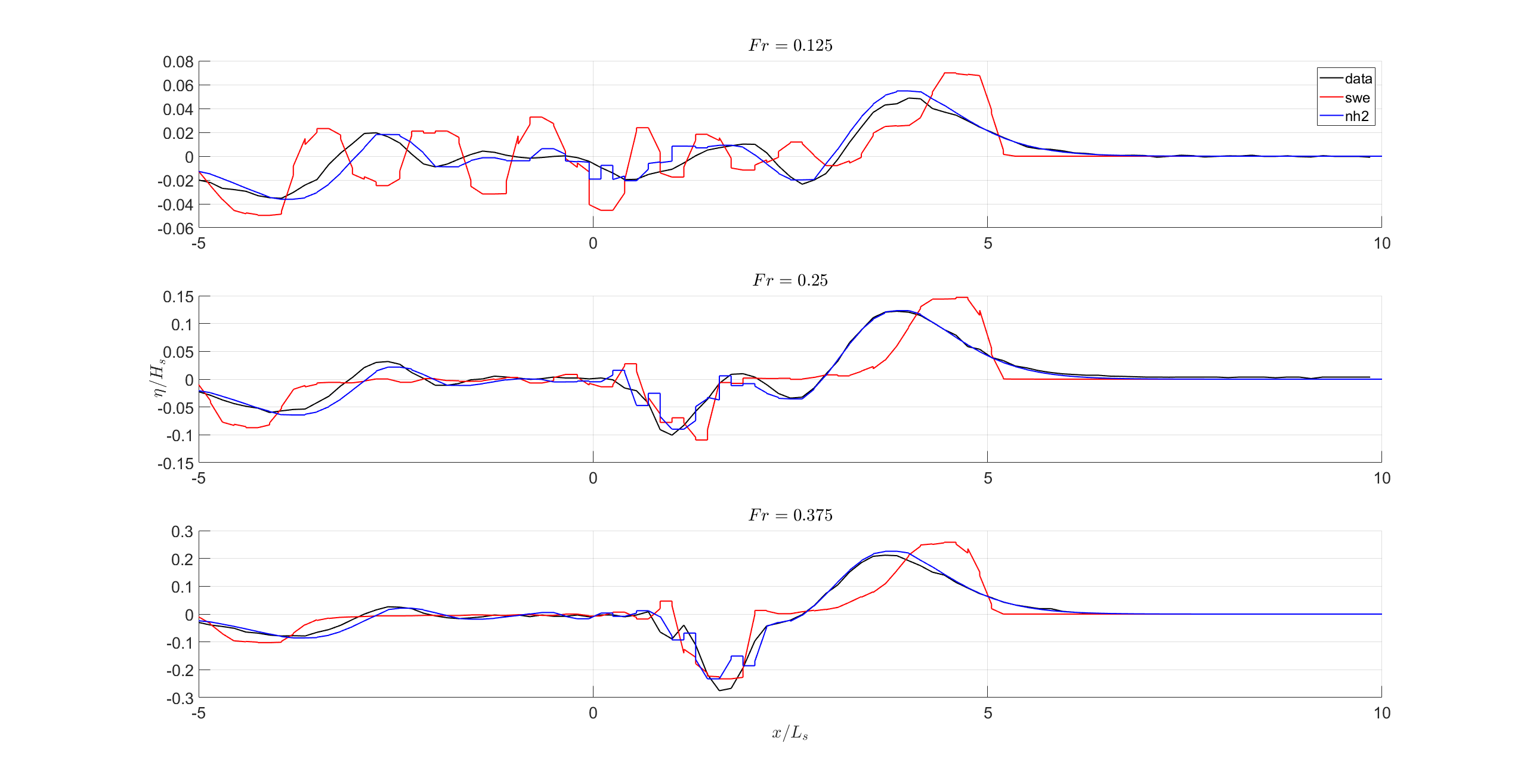}
    \caption{Comparison of the measured laboratory data (black) elevation with the simulation using hydrostatic SWE (red) and simulation using non-hydrostatic SWE (blue) for a sliding semi-elliptic plate for $Fr=0.125$ (top), $0.25$ (middle), and $0.375$ (bottom).}
    \label{fig:whittaker}
\end{figure}

For our simulation, following \cite{Jing2020}, the topography can be described as
\begin{equation}
    d(x,t) = h_0-H_s\left(1-\left(\frac{2(x-S(t))}{L_s}\right)^4\right), \quad S(t)-\frac{L_s}{2}<x<S(t)+\frac{L_s}{2},
\end{equation}
where $L_s = 0.5~m$ and $H_s = 0.026~m$ define the slide height and thickness respectively, while 
\begin{equation}
    S(t) = \begin{cases}
        \frac{1}{2}a_0t^2 & ,0\leq t\leq t_1,\\
        \frac{1}{2}a_0t_1^2+u_t(t-t_1) & ,t_1< t\leq t_2,\\
        \frac{1}{2}a_0t_1^2+u_t(t-t_1)-\frac{1}{2}a_0(t-t_2)^2 & ,t_2< t\leq t_3,\\
        \frac{1}{2}a_0t_1^2+u_t(t_3-t_1)-\frac{1}{2}a_0(t_3-t_2)^2 & ,t>t_3
    \end{cases}
\end{equation}
controls the slide movement.
Our computation will be compared to the data that was measured at $t = \frac{8}{\sqrt{g/L_s}}~s$. We use $\Delta t = 0.005 ~s$ and $\Delta x = 0.075~m$ for this simulation. The comparison with the measured data can be seen in Fig \ref{fig:whittaker}. In general, we can see how our non-hydrostatic model is more reliable than hydrostatic SWE, in terms of period and elevation.
For all these cases, our model successfully simulates the smooth leading waveform, while the hydrostatic SWE gives a relatively sharp form instead. Along with this shape, the hydrostatic SWE also gives a faster leading wave. In the lowest Froude number case ($Fr=0.125$), a significant difference can be observed, where the following wave gives an unexpected oscillatory behaviour in hydrostatic SWE. The non-hydrostatic model, on the other hand, captures a pretty accurate shape compared to the data, even though it overestimates the leading wave slightly. Moreover, for the higher Froude number cases ($Fr = 0.25$ and $0.375$) the leading wave amplitude is even better captured by our model. Meanwhile, the hydrostatic SWE fails to expose the dispersive effect on the following waves, where it gives a relatively flat surface instead. Hence, these comparisons may indicate the effectiveness of our non-hydrostatic model in simulating waves generated by accelerating landslides.

\subsection{Submerged landslide over a sloping beach}
Finally, we will implement our model to a phenomenon that happens over a slope, which is the case on most landslide tsunamis. This kind of test represents the effect of varying water depth on the onshore and offshore propagating waves. This case is based on a simulation as documented in \cite{lynett2002}, which is widely used as a benchmark for dispersive wave modelling \cite{FUHRMAN2009747,ashtiani2007, lynett2004}. Following \cite{lynett2002}, the bottom movement is given as 
\begin{equation}
    h(x,t)=x\tan\theta-\frac{\Delta h}{4}[1+\tanh(2\cos\theta(x-x_l(t)))][1-\tanh(2\cos\theta(x-x_r(t)))],
\end{equation}
\begin{equation}
    x_l(t) = x_c(t)-\frac{b}{2}\cos\theta,~x_r(t) = x_c(t)+\frac{b}{2}\cos\theta,
\end{equation}
\begin{equation}
    x_c(t) = x_0+S(t)\cos\theta,
\end{equation}
whereas $S(t) =  S_0\ln\left(\cosh\frac{t}{t_0}\right)$. This equation is slightly changed from the original study as suggested by \cite{FUHRMAN2009747}. We use the same coefficients as given in \cite{lynett2002}, which are $b = 1~m$, $\Delta h = 0.05~m$, $\theta = 6^{\circ}$, $S_0 = 4.712~m$ and $t_0 = 3.713~s$. Moreover, the centre of the sliding mass is located at $x_0=2.379~m$.

\begin{figure}[h]
    \centering
    \includegraphics[width=1 \linewidth]{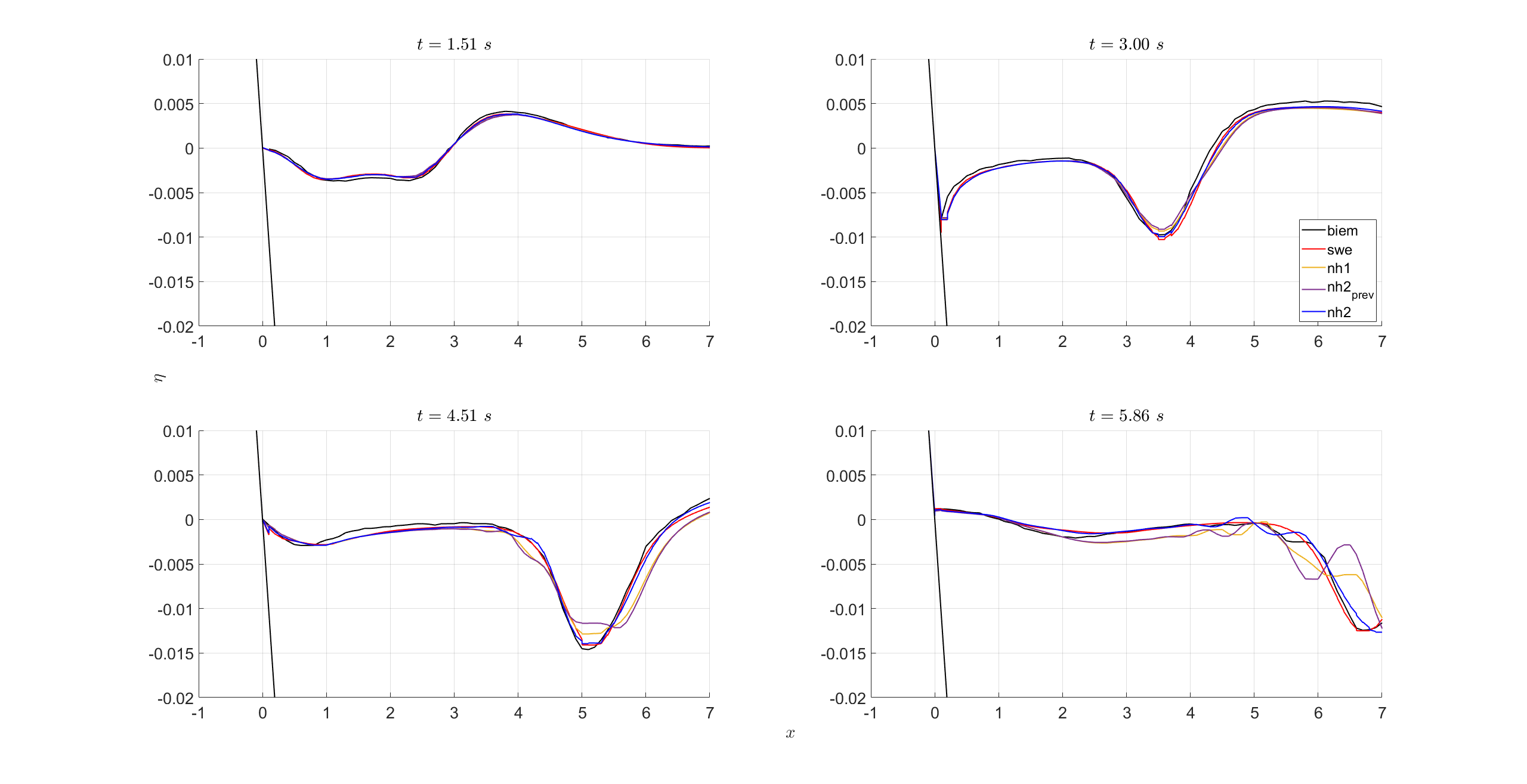}
    \caption{Comparison of the BIEM simulation (black) elevation with our simulation using hydrostatic SWE (red), linear non-hydrostatic SWE (yellow), simplified quadratic non-hydrostatic SWE (purple), and simulation using non-hydrostatic SWE (blue) for a sliding mass over sloping beach.}
    \label{fig:lynett}
\end{figure}

For our simulation, we use $\Delta t = 0.005~s$ and $\Delta x = 0.1~m$. The comparison is made at $t = 1.51, 3.00, 4.51$ and $5.86~s$ as depicted in Fig. \ref{fig:lynett}, where we compared against a result obtained with nonlinear and dispersive boundary integral equation method (BIEM) model, which is based on potential flow theory. The accuracy of this model is described in \cite{grilli1994}. Additionally, we also compare our derived model with a non-hydrostatic model with a linear pressure relation as well as a quadratic relation, with some neglected terms as in Eq. \ref{eq:linear_pressure} and \ref{eq:quadratic_pressure_old}.
It can be observed that for $t=1.51$ and $3.00~s$, all the simulations show good agreement with the BIEM results. At $t=4.51~s$, on the other hand, the linear and simplified quadratic relation gives a relatively wider waveform and inaccurate trough. Meanwhile, the hydrostatic SWE and our new non-hydrostatic model are well represented. Lastly, when the wave propagates to a deeper area at $t=5.86~s$, the linear and simplified non-hydrostatic models are even worsening as the waveform is not well captured. This test case shows the necessity to include the bottom derivative terms in the pressure relation.

The hydrostatic SWE are able to capture the waveform, yet they fail to simulate the dispersive effect that starts to occur in the BIEM simulation. Our new non-hydrostatic model shows a slightly faster phase speed while exposing a reasonable dispersive effect. Hence, our model here definitely outperforms the linear and simplified quadratic pressure relation, and it also reproduces the dispersive effect that occurs later stage. This test case demonstrates the new model's ability to solve such landslide-generated waves, with a more involved sloping bathymetry, which resembles more closely real landslide tsunami cases. 

\section{Conclusions}\label{sec5}
In this work, a non-hydrostatic SWE extension using a quadratic non-hydrostatic pressure relation for a moving bottom has been derived. Moreover, we proposed an alternative form of the vertical pressure profile, avoiding ambiguity due to a time derivative of the unknown velocity on the right-hand side of the non-hydrostatic equation. This allows for preserving the equivalence with the Boussinesq-type equations and solvability with the projection method.

The proposed set of non-hydrostatic equations requires the solution of an elliptic problem in each time step. However, with the re-formulation of the non-hydrostatic pressure relation the requirements for a flux formulation as given in \cite{Jeschke2017} are satisfied, allowing us to use the same procedure employing a local discontinuous Galerkin solution method. We demonstrate that the adopted flux does not destroy the representation of bottom topography and movement.

In order to validate our model, we have compared our numerical results with measured data and a high-fidelity simulation. First, we tested our model's ability to simulate waves generated by vertical thrust, which imitates a submerged earthquake-type phenomenon, or submarine volcano eruption, or caldera collapse. We also compared our results with a sliding obstacle on a flat or sloping bottom, which mimics a submerged landslide generation of waves. The first two benchmarks show that including non-hydrostatic pressure terms is crucial to improve hydrostatic SWE solutions in terms of waveform and period. We also showed that a quadratic pressure relation including the full bottom gradient terms is necessary to solve a more involved sloping bottom case accurately. Comparison with a linear and a simplified quadratic pressure relation formulation demonstrates the accuracy of the new approach.

We plan to develop this model formulation into a two-dimensional and more realistic simulation tool. In particular, we are interested in applying it to landslide tsunami scenarios with realistic bathymetry and source characteristics. Further developments should also consider friction and viscosity terms to overcome unrealistic wave amplitudes as exposed in the test cases with vertical thrust.

\section*{Acknowledgements}
The authors acknowledge the support of the Deutsche Forschungsgemeinschaft (DFG) within the Research Training Group GRK 2583 ”Modeling, Simulation and Optimization of Fluid Dynamic Applications”.

\bibliography{WileyNJD-AMA}
\end{document}